\documentclass[12pt]{article}
\usepackage{amssymb,amsfonts,amsmath, psfrag}
\usepackage{graphicx, colordvi, cite}
\usepackage{mathrsfs}
\usepackage{rotating}
\usepackage{hyperref}
\newtheorem{theo}{Theorem}

\newtheorem{coro}[theo]{Corollary}

\newtheorem{lem}[theo]{Lemma}

\makeatletter \@addtoreset{equation}{section}
\@addtoreset{theo}{section} \makeatother

\def\qed{\hfill \rule{4pt}{7pt}}
\def\pf{\noindent {\it Proof.} }
\begin{document}
\begin{center}
{\large\bf Asymptotic $r$-log-convexity and P-recursive sequences}
\end{center}

\begin{center}
Qing-Hu Hou$^{1}$ , Zuo-Ru Zhang$^{2}$\\[8pt]

Center for Applied Mathematics\\
Tianjin University\\
 Tianjin 300072, P. R. China\\[6pt]

$^{1}${\tt qh\_hou@tju.edu.cn}, $^{2}${\tt zhangzuoru@tju.edu.cn}
\end{center}
\vspace{0.3cm} \noindent{\bf Abstract.}
A sequence $\{ a_n \}_{n \ge 0}$ is said to be
asymptotically $r$-log-convex if it is $r$-log-convex for $n$ sufficiently large. We present a criterion on the asymptotical $r$-log-convexity based on the asymptotic behavior of $a_n a_{n+2}/a_{n+1}^2$.  As an application, we show that most P-recursive sequences are asymptotic $r$-log-convexity for any integer $r$ once they are log-convex.  Moreover, for a concrete integer $r$, we
present a systematic method to find the explicit integer $N$ such that a P-recursive sequence $\{a_n\}_{n \ge N}$ is $r$-log-convex. This enable us to prove the $r$-log-convexity of some combinatorial sequences.

\noindent {\bf Keywords:} asymptotic log-convexity, asymptotic $r$-log-convexity, P-recursive sequence

\section{Introduction}

Recall that a sequence $\{a_n\}_{n\geq 0}$ is said to be \emph{log-convex} (\emph{log-concave}, respectively)  if for all $n \ge 0$,
\begin{equation} \label{log-con}
a_{n} a_{n+2} - a_{n+1}^2 \ge 0 \ ( \le 0, respectively).
\end{equation}
Liu and Wang \cite{LW07} showed that many combinatorial sequences are log-convex. While for some sequences, the inequality \eqref{log-con} holds except for finitely many initial terms. For example, DeSalvo and Pak \cite{DeS15} proved that the partition function is log-concave for $n \ge 26$. This fact motives us to consider the \emph{asymptotic log-convexity} which requires \eqref{log-con} holding for $n$ sufficiently large. Recall that a sequence $\{a_n\}_{n\geq 0}$ is said to be \emph{$r$-log-convex} if
\begin{equation}\label{Lk}
 \{ \mathscr{L} a_n \}_{n \ge 0}, \  \{ \mathscr{L}^2 a_n \}_{n \ge 0}, \ \ldots,\   \{ \mathscr{L}^r a_n \}_{n \ge 0}
\end{equation}
are all non-negative sequences, where
\[
 \mathscr{L} a_n = a_{n} a_{n+2} - a_{n+1}^2 \quad \mbox{and} \quad \mathscr{L}^{\, k} a_n = \mathscr{L} (\mathscr{L}^{\,k-1} a_n).
\]
We also define the \emph{asymptotic $r$-log-convexity} by requiring the sequences in \eqref{Lk} to be non-negative for $n$ sufficiently large.

It is straightforward to see that a positive sequence $\{a_n\}_{n\geq 0}$ is asymptotic $r$-log-convexity if the ratio $a_{n} a_{n+2}/ a_{n+1}^2$ is asymptotic to $1+c/n^\alpha$
for certain positive numbers $c$ and $\alpha$. By investigating the relation between
\[
\frac{a_n a_{n+2}}{a_{n+1}^2} \quad \mbox{and} \quad
\frac{\mathscr{L}a_n \mathscr{L}a_{n+2}}{(\mathscr{L}a_{n+1})^2},
\]
we extend this observation to the case of  asymptotic $r$-log-convexity and present a criterion.

With this criterion, we studied the asymptotic $r$-log-convexity of the \emph{P-recursive} sequences. Recall that a P-recursive sequence of order $d$ satisfies a recurrence relation of the form
\[
a_n = r_1(n)a_{n-1} + r_2(n)a_{n-2} + \cdots + r_d(n) a_{n-d},
\]
where $r_i(n)$ are rational functions of $n$ (see \cite[Section 6.4]{Sta99}).
Do\v{s}li\'{c} and Veljan \cite{Dos08} presented a method on proving the log-convexity of the P-recursive sequences. Chen and Xia \cite{Chen11} gave a criterion for the $2$-log-convexity of the P-recursive sequences of order $2$. Zhu \cite{Zhu13} considered the ratio log-convexity of the P-recursive sequences of order $2$. All of these methods depend on manually finding suitable bounds of the ratio $a_{n+1}/a_n$.

Our approach is based on the asymptotic behavior of the P-recursive sequences given by Birkhoff and Trjitzinsky \cite{Bir32} and developed by Wimp and Zeilberger \cite{WZ85}.
They showed that a P-recursive sequence $\{a_n\}_{n \ge 0}$ is asymptotically equal to a linear combination of terms of the form
\begin{equation}
e^{Q(\rho, n)} s(\rho, n),
\end{equation}
where
\[
Q(\rho, n) = \mu_0n\log n+ \sum_{j=1}^\rho \mu_j n^{j/\rho},
\]
\[
s(\rho, n) = n^r \sum_{j=0}^{t-1} (\log n)^j \sum_{s=0}^{M-1} b_{sj}n^{-s/\rho},
\]
with $\rho,t,M$ being positive integers and $\mu_j, r,  b_{sj}$ being complex numbers.
We show that when $t=1$, the asymptotic $r$-log-convexity of $\{ a_n\}_{n \ge 0}$ can be deduced by computing $Q(\rho,n)$ and $s(\rho,n)$.  As an example, we confirm the second part of Conjecture 5.4 in \cite{Chen11}.

The asymptotic expression of $a_n$ also provides good bounds for the ratio  $a_{n+1}/a_n$. As we have mentioned, these bounds play important roles in applying the methods in \cite{Dos08, Chen11, Zhu13}. Based on the asymptotic expression, we present a systematic method on proving the $r$-log-convexity of the P-recursive sequences. As an application, we show that the number of $n \times n$ $(0,1)$ matrices with row sum and column sum $3$ is $2$-log-convex for $n \ge 8$.

This paper is organized as follows.  In Section 2, by studying the relation of $\mathscr{R}^2 a_n$ and $\mathscr{R}^2 \mathscr{L} a_n$, we give a criterion of the asymptotic $r$-log-convexity for sequences whose ratios $a_n a_{n+2}/a_{n+1}^2$ have the Puiseux-type approximations. We show that most P-recursive sequences are of this kind and thus confirm a conjecture on the asymptotic $r$-log-convexity of Motzkin numbers. In Section 3, we give a  method to prove the $r$-log-convexity of  P-recursive sequences by finding the explicit $N$ such that  $\{a_n\}_{n\geq N}$ is $r$-log-convex.

\section{Asymptotic $r$-log-convexity}

In this section, we will consider a kind of positive sequences $\{a_n\}_{n \ge 0}$  such that the ratio $a_n a_{n+2}/a_{n+1}^2$ is asymptotic to
\[
1 + \frac{c_1}{n^{\alpha_1}} +  \frac{c_2}{n^{\alpha_2}} + \cdots +  \frac{c_m}{n^{\alpha_m}},
\]
where $0 < \alpha_1 < \alpha_2 < \cdots < \alpha_m$ and $c_i$ are real numbers. We will give a criterion for the asymptotic $r$-log-convexity of these sequences. Moreover, we will show that most P-recursive sequences satisfy this condition.

Let $\mathscr{R} a_n = a_{n+1}/a_n$.
We see that the log-convexity of a positive sequence $\{a_n\}_{n \ge 0}$ is equivalent to the statement that
\[
\mathscr{R}^2 a_n = \frac{a_n a_{n+2}}{a_{n+1}^2} \ge 1.
\]
To study the $r$-log-convexity, we firstly investigate the relation between
$\mathscr{R}^2 a_n$ and $\mathscr{R}^2 \mathscr{L} a_n$.

\begin{lem}\label{re}
Denote $s_n = \mathscr{R}^2 a_n$. Then
\[
\mathscr{R}^2 \mathscr{L} a_n = s_{n+1}^2 \frac{(s_n-1)(s_{n+2}-1)}{(s_{n+1}-1)^2}.
\]
\end{lem}
\pf By direct computation, we have
\begin{align*}
 {\mathscr R}^2{\mathscr L}{a_n} & = \frac{(a_{n+3}a_{n+1}-a_{n+2}^2)(a_{n+1}a_{n-1}-a_{n}^2)}{(a_{n+2}a_n-a_{n+1}^2)^2}\\
                           & = \frac{{a_{n+2}^2}a_{n}^2}{{a_{n+1}^2}} \frac{(\frac{a_{n+3}a_{n+1}}{a_{n+2}^2}-1)(\frac{a_{n+1}a_{n-1}}{a_n^2}-1)}{(\frac{a_{n+2}a_n}{a_{n+1}^2}-1)^2}\\
                           & = {s_{n+1}^2}\frac{(s_{n+2}-1)(s_n-1)}{(s_{n+1}-1)^2},
\end{align*}
completing the proof. \qed

Let $\{f_n\}_{n\geq0}$ be a sequence of real numbers. Suppose that there exist
real numbers $c_i, \alpha_i$ with
\[
\alpha_0 < \alpha_1 < \cdots < \alpha_m
\]
such that
\[
\lim_{n \to \infty} n^{\alpha_m} \left( f_n -  \sum_{i=0}^m \frac{c_i}{n^{\alpha_i}} \right) = 0.
\]
We call
\[
g_n = \sum_{i=0}^m \frac{c_i}{n^{\alpha_i}}
\]
a \emph{Puiseux-type approximation} of $f_n$ and denote $f_n \approx g_n$.
We will see that the first few terms of $g_n$ and the exponent $\alpha_m$ play key roles in our discussion. Therefore, we will abbreviate the trailing terms and use the standard little-o notation to write
\[
f_n = \frac{c_0}{n^{\alpha_0}} + \frac{c_1}{n^{\alpha_1}} + \cdots + o\left( \frac{1}{n^{\alpha_m}} \right).
\]

Based on Lemma~\ref{re}, we can derive a criterion for the asymptotic $r$-log-convexity in terms of the Puiseux-type approximation of $\mathscr{R}^2 a_n$.

\begin{theo} \label{asy}
Let $\{a_n\}_{n\geq0}$ be a positive sequence such that $\mathscr{R}^2 a_n$ has a Puiseux-type approximation
\begin{equation}\label{sn}
\mathscr{R}^2 a_n = 1 + \frac{c}{n^{\alpha}} + \cdots + o \left( \frac{1}{n^{\beta}} \right),
\end{equation}
where $c>0$ and $0<\alpha \le \beta$. Then $\{a_n\}_{n\geq0}$ is asymptotically $r$-log-convex for
\begin{equation}\label{r-log}
r= \begin{cases}
	\left\lfloor \beta/\alpha \right\rfloor, & \mbox{if $\alpha < 2$}, \\[5pt]
	\left\lfloor \frac{\beta-\alpha}{2} \right\rfloor+1, & \mbox{if $\alpha \ge 2$},
    \end{cases}
\end{equation}
where $\lfloor x \rfloor$ denotes the maximal integer less than or equal to $x$.
\end{theo}
\pf
Let $s_n = \mathscr{R}^2 a_n$.  By \eqref{sn} we immediately derive that
$s_n \ge 1$ for $n$ sufficiently large and hence $\{a_n\}$ is asymptotically log-convex.
We will compute a Puiseux-type approximation of $\mathscr{R}^2 \mathscr{L} a_n$ by aid of Lemma~\ref{re}.

Noting that
\[
\frac{1}{(n+\gamma)^\alpha} = \frac{1}{n^\alpha}  + \cdots + o \left( \frac{1}{n^{\beta}} \right),
\]
we have
\[
s_{n+1}^2  =  \left( 1 + \frac{c}{n^{\alpha}} + \cdots + o \left( \frac{1}{n^{\beta}} \right) \right)^2
 =
1 + \frac{2c}{n^{\alpha}} + \cdots  + o \left( \frac{1}{n^{\beta}} \right).
\]

Denote
\[
f_n = \frac{n^{\alpha}}{c} (s_n-1).
\]
We have
\begin{equation}
f_n = 1 + \frac{d}{n^{\gamma}} + \cdots + o \left( \frac{1}{n^{\beta-\alpha}} \right),
\end{equation}
for certain real numbers $d$ and $\gamma$. To estimate the ratio $f_n f_{n+2}/f_{n+1}^2$, we consider its logarithm. By the taylor expansion of logarithm function, we see that
\[
\log f_n = \frac{d}{n^{\gamma}} + \cdots + o \left( \frac{1}{n^{\beta-\alpha}} \right).
\]
Notice that for any $\theta>0$,
\begin{align*}
& \frac{1}{n^\theta} + \frac{1}{(n+2)^\theta} - \frac{2}{(n+1)^\theta} \\[5pt]
& = \frac{1}{n^\theta}  + \left( \frac{1}{n^\theta} - \frac{2 \theta}{n^{\theta+1}}
+ \frac{2\theta(\theta+1)}{n^{\theta+2}} \right)  - 2  \left( \frac{1}{n^\theta} - \frac{\theta}{n^{\theta+1}}
+ \frac{\theta(\theta+1)}{2 n^{\theta+2}} \right) \\
& \qquad + \cdots + o \left( \frac{1}{n^{\beta-\alpha}} \right) \\
& = \frac {\theta(\theta+1)}{n^{\theta+2}} + \cdots + o \left( \frac{1}{n^{\beta-\alpha}} \right).
\end{align*}
We thus derive that
\[
\log \frac{f_n f_{n+2}} {f_{n+1}^2} = \frac{d \gamma (\gamma+1)}{n^{\gamma+2}} + \cdots + o \left( \frac{1}{n^{\beta-\alpha}} \right),
\]
and hence
\begin{equation}\label{r2}
 \frac{f_n f_{n+2}} {f_{n+1}^2} = 1 +  \frac{d \gamma (\gamma+1)}{n^{\gamma+2}} + \cdots + o \left( \frac{1}{n^{\beta-\alpha}} \right).
\end{equation}
Therefore,
\begin{align*}
& \frac{(s_n-1)(s_{n+2}-1)}{(s_{n+1}-1)^2}   = \frac{(n+1)^{2\alpha}} {n^{\alpha}(n+2)^{\alpha} } \frac{g_n g_{n+2} }{g_{n+1}^2} \\
& = \left( 1 + \frac{\alpha}{n^2} + \cdots + o \left( \frac{1}{n^{\beta}} \right) \right) \cdot \left( 1 +  \frac{d \gamma (\gamma+1)}{n^{\gamma+2}} + \cdots + o \left( \frac{1}{n^{\beta-\alpha}} \right)  \right) \\
& = 1 + \frac{\alpha}{n^2} + \cdots + o \left( \frac{1}{n^{\beta-\alpha}} \right).
\end{align*}

Now by Lemma~\ref{re} we deduce that
\[
\mathscr{R}^2 \mathscr{L} a_n
 = \left( 1 + \frac{2c}{n^{\alpha}} + \cdots  + o \left( \frac{1}{n^{\beta}} \right) \right) \left( 1 + \frac{\alpha}{n^2} + \cdots + o \left( \frac{1}{n^{\beta-\alpha}} \right) \right).
\]

We firstly consider the case of $\alpha<2$. If $\beta-\alpha \ge \alpha$, then we obtain a Puiseux-type approximation of $\mathscr{R}^2 \mathscr{L} a_n$:
\[
1 + \frac{2c}{n^{\alpha}} + \cdots  + o \left( \frac{1}{n^{\beta -\alpha}} \right).
\]
Since $c>0$, we derive that $\{a_n\}_{n \ge 0}$ is asymptotically $2$-log-convex.
Repeating the above discussion, we finally derive that $\{a_n\}_{n \ge 0}$ is asymptotically $r$-log-convex for $r \le \lfloor \beta/\alpha \rfloor$.

Then we consider the case of $\alpha=2$. In this case, if $\beta-\alpha \ge 2$, we will obtain a Puiseux-type approximation of $\mathscr{R}^2 \mathscr{L} a_n$:
\[
1 + \frac{(2c+2)}{n^2} + \cdots + o \left( \frac{1}{n^{\beta-2}} \right).
\]
Also, by iterating the above discussion, we obtain that $\{a_n\}_{n \ge 0}$ is asymptotically $r$-log-convex for $r \le \lfloor \beta/2 \rfloor$.

Finally we consider the case of $\alpha>2$. In this case, when $\beta-\alpha \ge 2$, we will obtain a Puiseux-type approximation of $\mathscr{R}^2 \mathscr{L} a_n$ as
\[
1 + \frac{\alpha}{n^2} + \cdots + o \left( \frac{1}{n^{\beta-\alpha}}  \right).
\]
By iterating the above discussion, we deduce that $\{a_n\}_{n \ge 0}$ is asymptotically $r$-log-convex for $r \le \lfloor (\beta-\alpha)/2 \rfloor+1$. Combining the results for $\alpha=2$ and $\alpha>2$, we arrive at \eqref{r-log}. \qed

\noindent {\it Remark.} We see that the asymptotic $r$-log-convexity depends only on the first term of the Puiseux-type approximation of $\mathscr{R}^2 a_n$. Suppose that for any real number $\beta$, $\mathscr{R}^2 a_n$ has a Puiseux-type approximation of form \eqref{sn}. Then the sequence $\{a_n\}_{n \ge 0}$ is asymptotically  $r$-log-convex for any integer $r$ once it is log-convex.

In order to apply Theorem~\ref{asy}, we need to show that $\mathscr{R}^2 a_n$ has a Puiseux-type approximation of the form \eqref{sn}. We will show that this is the case for most P-recursive sequences.

Recall that a formal solution to a polynomial recurrence relation is of the form
\[
f_n = e^{Q(\rho, n)} s(\rho, n),
\]
where
\[
Q(\rho, n) = \mu_0n\log n+ \sum_{j=1}^\rho \mu_jn^{j/\rho},
\]
\[
s(\rho, n) = n^r \sum_{j=0}^{t-1} (\log n)^j \sum_{s=0}^{\infty} b_{sj}n^{-s/\rho},
\]
with $\rho,t$ being positive integers and $\mu_j, r,  b_{sj}$ being complex numbers. We will see that the ratio $\mathscr{R}^2 f_n$ has a good asymptotic behavior when $t=1$.

\begin{theo}\label{P-rec}
Suppose that $\{a_n\}_{n \ge 0}$ is a P-recursive sequence whose asymptotic expression is
\[
a_n = e^{Q(\rho, n)} \cdot n^r \left( \sum_{s=0}^{M} b_{s} n^{-s/\rho} + o \left( n^{-M/\rho} \right) \right),
\]
where
\[
Q(\rho, n) = \mu_0n\log n+ \sum_{j=1}^\rho \mu_j n^{j/\rho},
\]
with $\rho,M$ being positive integers, $\mu_j, r,  b_{s}$ being real numbers and $b_0 \not= 0$.
Then there exist real numbers $c_1,\ldots,c_M$ such that
\begin{equation}\label{R2a}
\mathscr{R}^2 a_n \approx 1 + \sum_{i=1}^{M} \frac{c_i}{n^{i/\rho}}.
\end{equation}
\end{theo}

\pf Noting that the operator $\mathscr{R}$ is multiplicative, we may consider each factor of $a_n$ separetively.

We firstly consider the factor $f_n = e^{\mu_0n\log n}$. We have
\begin{align*}
\log (\mathscr{R}^2 f_n) & = \mu_0 \big( (n+2)\log(n+2)+n \log n - 2(n+1)\log(n+1) \big) \\
&= \mu_0 \left( (n+2)\log\left( 1+\frac{2}{n} \right) - 2(n+1)\log\left( 1+\frac{1}{n} \right) \right) \\
&  = \mu_0 \left( \frac{1}{n} - \frac{1}{n^2} + \cdots + o\left( \frac{1}{n^M} \right) \right).
\end{align*}
Therefore,
\[
\mathscr{R}^2 f_n = 1 + \frac{\mu_0}{n} - \frac{\mu_0^2/2-\mu_0}{n^2} + \cdots + o\left( \frac{1}{n^M} \right) .
\]

Then we consider the factor $f_n = e^{\mu_j n^{j/\rho}}$. We have
\begin{align*}
\log (\mathscr{R}^2 f_n) & = \mu_j \big( (n+2)^{j/\rho} + n^{j/\rho} - 2(n+1)^{j/\rho} \big) \\
&=  \mu_j n^{j/\rho}  \left( \left( 1+\frac{2}{n} \right)^{j/\rho} + 1 - 2  \left( 1+\frac{1}{n} \right)^{j/\rho} \right) \\
&  = \mu_j n^{j/\rho} \left( \frac{(j/\rho)(j/\rho-1)}{n^2}  + \cdots + o\left( \frac{1}{n^{M+1}} \right) \right).
\end{align*}
Therefore,
\[
\mathscr{R}^2 f_n = 1 + \frac{\mu_j (j/\rho)(j/\rho-1)}{n^{2-j/\rho}} + \cdots + o\left( \frac{1}{n^{M}} \right) .
\]

Next we consider the factor $f_n=n^r$. We have
\[
\mathscr{R}^2 f_n = \left( 1+\frac{2}{n} \right)^r \left(1+\frac{1}{n} \right)^{-2r} = 1 - \frac{r}{n^2} + \cdots + o\left( \frac{1}{n^{M}}  \right).
\]

Finally, we consider the factor
\[
f_n = \sum_{s=0}^{M} b_{s} n^{-s/\rho} + o \left( n^{-M/\rho} \right).
\]
The ratio $f_n f_{n+2}/f_{n+1}^2$ has been estimated in the proof of Theorem~\ref{asy} for $b_0=1$. By dividing $f_n$ by $b_0$, we get
\[
\mathscr{R}^2 f_n = 1 + \frac{(b_1/b_0) \cdot (1/\rho) \cdot (1/\rho+1)}{n^{1/\rho+2}} + \cdots+ o \left( n^{-M/\rho} \right).
\]

Combining the factors together, we finally arrive at Equation~\eqref{R2a}. \qed

Now we give an example to illustrate how to apply Theorems \ref{asy} and \ref{P-rec}.

\begin{coro}
Let $M_n$ be the $n$-th Motzkin number. Then the sequence $\{M_n\}_{n \ge 0}$ is asymptotically $r$-log-convex for any positive integer $r$.
\end{coro}
\pf
It is known that the Motzkin numbers satisfy the recurrence
\begin{equation}\label{rec-M}
(n+2)M_n=(2n+1)M_{n-1}+(3n-3)M_{n-2},
\end{equation}
where $n\geq2$ and $M_0=M_1=1$, see Aigner\cite{Aig98}.

We can get the asymptotic expression of $M_n$ from the recurrence relation by considering the asymptotic behavior of the ratio $M_{n+1}/M_n$. Zeilberger has implemented a {\tt Maple} package {\tt AsyRec} \cite{Zeil16}. We have our own {\tt Mathematica} implementation which is accessible at {\tt http://cam.tju.edu.cn/\~\,hou/preprints.html}.

By these packages, we deduce two formal solutions to the recursion \eqref{rec-M}:
\[
3^n n^{-3/2} \sum_{i=0}^\infty a_i n^{-i} \quad \mbox{and} \quad (-1)^n n^{-3/2} \sum_{i=0}^\infty b_i n^{-i},
\]
where $a_i,b_i$ are real numbers. Since the second solution tends to zero when $n$ tends to infinite, we thus derive that $M_n$ has a Puiseux-type approximation of the form
\[
M_n \approx 3^n n^{-3/2} \sum_{i=0}^K a_i n^{-i},
\]
for any positive integer $K$. By comparing the first few terms, we compute that
\[
M_n \approx C \cdot \frac{3^n}{n^{3/2}} \left(1- \frac{39}{16 n} \right),
\]
for some constant $C$. Hence,
\[
\mathscr{R}^2 M_n \approx 1 + \frac{3}{2n^2}.
\]
We thus derive that $\{M_n\}_{n \ge 0}$ is asymptotically $r$-log-convex for any positive integer $r$. \qed

We remark that $\{M_n\}_{n \ge 0}$ is not $2$-log-convex but $\{M_n \}_{n \ge 6}$ is.

With the same discussion, besides the log-behavior of Catalan-Larcombe-French sequence given by Sun and Zhao \cite{Sun16}, we find that the Catalan-Larcombe-French sequence, the Fine numbers and the Franel numbers of orders $3$--$6$ are asymptotically $r$-log-convex for any integer $r$. This confirms the second part of the conjecture posed by Chen and Xia \cite{Chen11}.

\section{A method on proving the $r$-log-convexity}
In the previous section, we see how to obtain the asymptotic $r$-log-convexity of P-recursive sequences. For a concrete integer $r$, we can indeed prove the $r$-log-convexity by finding the explicit $N$ such that $\{a_n\}_{n\geq N}$ is $r$-log-convex.

Assume that $\{a_n\}_{n \ge 0}$ is a positive sequence which satisfies the recurrence relation
\begin{equation}\label{recR}
a_{n+t}=R_0(n)a_n+R_1(n)a_{n+1}+ \cdots +R_{t+1}(n)a_{n+t-1}.
\end{equation}
Suppose that for any positive integer $m$, there exists a Puiseux-type approximation
\[
a_{n+1}/a_n  \approx \sum_{i=m_0}^m \frac{c_i}{n^{i/\rho}},
\]
where $\rho$ is a positive integer, $m_0$ is an integer and $c_i$ are real numbers.
We say $\{a_n\}_{n \ge 0}$ is a \emph{bound preserving} sequence if for $m$ sufficiently large and
\[
 f_n = \sum_{i=m_0}^m \frac{c_i}{n^{i/\rho}} - \frac{1}{n^{m/\rho}}, \quad
 g_n = \sum_{i=m_0}^m \frac{c_i}{n^{i/\rho}} + \frac{1}{n^{m/\rho}},
\]
we can compute an integer $N_0$ such that
\[
\frac{R_0(n)}{u^{(0)}_n u^{(0)}_{n+1} \cdots u^{(0)}_{n+t-2}}+
\frac{R_1(n)}{u^{(1)}_{n+1} \cdots u^{(1)}_{n+t-2}} + \cdots+R_{t-1}(n) \geq f_{n+t-1}^{(0)}, \quad \forall\, n \ge N_0,
\]
and
\[
\frac{R_0(n)}{v^{(0)}_n v^{(0)}_{n+1} \cdots v^{(0)}_{n+t-2}}+
\frac{R_1(n)}{v^{(1)}_{n+1} \cdots v^{(1)}_{n+t-2}} + \cdots+R_{t-1}(n) \leq g_{n+t-1}^{(0)}, \quad \forall\, n \ge N_0,
\]
where
\[
u^{(i)}_n=
\begin{cases}
g_n, & \mbox{if $R_i(n)>0$ for $n$ sufficiently large,}\\[5pt]
f_n, &  \mbox{if $R_i(n)<0$ for $n$ sufficiently large.}
\end{cases}
\]
and
\[
v^{(i)}_n=\begin{cases}
f_n, & \mbox{if $R_i(n)>0$ for $n$ sufficiently large,}\\[5pt]
g_n, & \mbox{if $R_i(n)<0$ for $n$ sufficiently large.}
\end{cases}
\]

 For a bound preserving sequence $\{a_n\}_{n \ge 0}$, we can find the explicit $N$ such that  $\{a_n \}_{n \ge N}$ is $r$-log-convex, when it is asymptotically log-convex.

\begin{theo}
Let $\{a_n\}_{n \ge 0}$ be a bound preserving sequence. Assume that
\[
\mathscr{R}^2 a_n \approx 1 + \frac{c}{n^\alpha},
\]
where $c$ and $\alpha$ are positive numbers.
Then for each positive integer $r$, we can compute an integer $N$ such that $\{a_n\}_{n \ge N}$ is $r$-log-convex.
\end{theo}
\pf
Our idea is to find good bounds of $\mathscr{R}^2 \mathscr{L}^i a_n$ for $i=0, 1, \ldots, r-1$ based on the Puiseux-type approximation of $\mathscr{R}^2 a_n$.

Denote
\[
s_n^{(i+1)} = \mathscr{R}^2 \mathscr{L}^i a_n, \quad i=0,1,\ldots,r-1.
\]
By the hypotheses, for any integer $m$ and $i$, there exists a Puiseux-type approximation
\[
s_n^{(i)} \approx 1 + \sum_{j=1}^m \frac{c_j^{(i)}}{n^{j/\rho}}.
\]

Since $s_n^{(1)} \approx 1+\frac{c}{n^{\alpha}}$ with $c>0$ and $\alpha>0$, we know that $s_n^{(i)} \approx 1 +\frac{c_i}{n^{\alpha_i}}$ for some $c_i >0$ and $\alpha_i >0$. We firstly find bounds of $s_n^{(r-1)}$ to ensure that $s_n^{(r)} \ge 1$. To this end, we search for $k_1$ and $N_1$ such that
\begin{equation}\label{f1}
f_n^{(r-1)} =1+\sum_{i=1}^{k_1}\frac{c_i^{(r-1)}}{n^{i/\rho}}-\frac{1}{n^{k_1/\rho}},
\end{equation}
\begin{equation}\label{g1}
g_n^{(r-1)} = 1+\sum_{i=1}^{k_1}\frac{c_i^{(r-1)}}{n^{i/\rho} }+ \frac{1}{n^{k_1/\rho}},
\end{equation}
and
\begin{equation}\label{snr}
\frac{(f_n^{(r-1)}-1)(f_{n+2}^{(r-1)}-1)}{(1-1/g_{n+1}^{(r-1)})^2}-1\geq 0,\quad \forall\,  n\geq N_1.
\end{equation}
Since $s_n^{(r)} \approx 1 + \frac{c_r}{n^{\alpha_r}}$ for some $c_r >0$ and $\alpha_r>0$,
there always exist $k_1$ and $N_1$ such that \eqref{f1}--\eqref{snr} holds. On the contrary, once we have
\begin{equation}\label{bounds}
f_n^{(r-1)} \le s_n^{(r-1)} \le g_n^{(r-1)}, \quad \forall \, n \ge N_1,
\end{equation}
we will deduce that
\[
s_n^{(0)} \ge \frac{(f_n^{(r-1)}-1)(f_{n+2}^{(r-1)}-1)}{(1-1/g_{n+1}^{(r-1)})^2} \ge 1, \quad \forall \, n \ge N_1.
\]

Next we shall find bounds of $s_n^{(r-2)}$ to ensure that \eqref{bounds} holds.
Similar to the above discussion, we search for $k_2$ and $N_2$ such that
\[
f_n^{(r-2)} =1+\sum_{i=1}^{k_2}\frac{c_i^{(r-2)}}{n^{i/\rho}}-\frac{1}{n^{k_2/\rho}},
\]
\[
g_n^{(r-2)} = 1+\sum_{i=1}^{k_2}\frac{c_i^{(r-2)}}{n^{i/\rho} }+ \frac{1}{n^{k_2/\rho}},
\]
and
\[
\frac{(f_{n+2}^{(r-2)}-1)(f_n^{(r-2)}-1)}{(1-1/g_{n+1}^{(r-2)})^2}\geq f_n^{(r-1)},\quad\forall\, n\geq N_2,
\]
\[
\frac{(g_{n+2}^{(r-2)}-1)(g_n^{(r-2)}-1)}{(1-1/f_{n+1}^{(r-2)})^2}\leq g_n^{(r-1)},\quad \forall\, n\geq N_2.
\]

Repeat the above  process until we find lower and upper bounds of $s_n^{(1)}$:
\[
f_n^{(1)} = 1+\sum_{i=1}^{k_{r-1}}\frac{c_i^{(1)}}{n^{i/\rho}}-\frac{1}{n^{k_{r-1}/\rho}},
\]
and
\[
g_n^{(1)} = 1+\sum_{i=1}^{k_{r-1}}\frac{c_i^{(1)}}{n^{i/\rho}}+\frac{1}{n^{k_{r-1}/\rho}}.
\]

Now we will find bounds of $a_{n+1}/a_n$. We search for $k_r$ and $N_r$ such that
\begin{equation}
f_n^{(0)}=\sum_{i=m_0}^{k_r}\frac{c_i}{n^{i/\rho}}-\frac{1}{n^{k_r/\rho}},
\end{equation}
\begin{equation}
g_n^{(0)}=\sum_{i=m_0}^{k_r}\frac{c_i}{n^{i/\rho}}+\frac{1}{n^{k_r/\rho}},
\end{equation}
and
\begin{equation}
f_{n+1}^{(0)}/g_n^{(0)} \geq f_n^{(1)}, \quad \forall\, n\geq N_r,
\end{equation}
\begin{equation}
g_{n+1}^{(0)} /f_n^{(0)} \leq g_n^{(1)}, \quad \forall\, n\geq N_r.
\end{equation}

Finally, we determine $N_0$ such that
\[
f_n^{(0)} \le a_{n+1}/a_n \le g_n^{(0)}, \quad \forall\, n \ge N_0,
\]
by aid of the recurrence relation \eqref{recR}. For this purpose, we search for $N_0$ such that

$1)$ $R_0(n), \ldots, R_{t-1}(n)$ are uniformly negative or positive for $n\geq N_0.$

$2)$ For $n\geq N_0$, we have
\begin{equation}\label{u}
\frac{R_0(n)}{u^{(0)}_n u^{(0)}_{n+1} \cdots u^{(0)}_{n+t-2}}+
\frac{R_1(n)}{u^{(1)}_{n+1} \cdots u^{(1)}_{n+t-2}} + \cdots+R_{t-1}(n) \geq f_{n+t-1}^{(0)},
\end{equation}
where
\[
u^{(i)}_n=
\begin{cases}
g_n^{(0)}, & \mbox{if $R_i(n)>0$,}\\[5pt]
f_n^{(0)}, &  \mbox{if $R_i(n)<0$.}
\end{cases}
\]

$3)$ For $n\geq N_0$, we have
\begin{equation}\label{v}
\frac{R_0(n)}{v^{(0)}_n v^{(0)}_{n+1} \cdots v^{(0)}_{n+t-2}}+
\frac{R_1(n)}{v^{(1)}_{n+1} \cdots v^{(1)}_{n+t-2}} + \cdots+R_{t-1}(n) \leq g_{n+t-1}^{(0)},
\end{equation}
where
\[
v^{(i)}_n=\begin{cases}
f_n^{(0)} & \mbox{if $R_i(n)>0$,}\\[5pt]
g_n^{(0)} & \mbox{if $R_i(n)<0$.}
\end{cases}
\]
$4)$ For $n=N_0,N_0+1,...,N_0+t-1$, we have
\[
f_n^{(0)}\leq a_{n+1}/a_n \leq g_n^{(0)}.
\]
Since $\{a_n\}_{n \ge 0}$ is an order preserving sequence, we can find out $N_0$ explicitly.

Set $N=\max\{N_0, N_1,...,N_r\} $. Then the sequence $\{a_n\}_{n\geq N}$ is $r$-log-convex. \qed

We conclude by giving an example.

Let $H_n(3)$ denote the number of $n\times n$ $(0,1)$-matrices with row and column sum 3. By the theory of symmetric functions, one can derive that
\begin{multline*}
  H_n-\frac{9n^2+42n+37}{(n+1)^2(3n+10)} H_{n+1}-\frac{3(3n+13)}{(n+2)(n+1)^2(3n+10)} H_{n+2} \\
 -\frac{9(3n^2+19n+32)}{(n+3)(n+2)^2(n+1)^2(3n+10)} H_{n+3} \\
 +\frac{12(3n+7)}{(n+4)(3n+10)(n+3)^2(n+2)^2(n+1)^2} H_{n+4}=0 .
 \end{multline*}
 Yang conjectured that $\{H_n(3)\}_{n \ge 2}$ is log-convex (private communication). We have the following stronger result.

\begin{theo}
The sequence $\{H_n(3)\}_{n \ge 8}$ is $2$-log-convex.
\end{theo}
\pf
By our {\tt Mathematica} package, we find that
\[
H_n(3) = C \cdot e^{-3n}n^{3n+1/2} \left( \frac{3}{4} \right)^n \left( 1-\frac{55}{36n} -\frac{47}{2592 n^2} + \frac{1045741}{1399680 n^3} + o \left( \frac{1}{n^3} \right)\right),
\]
for some constant $C$. Therefore,
\[
r_n=\frac{H_{n+1}(3)}{H_n(3)}\approx \frac{3}{4}n^3+\frac{3}{2}n^2+\frac{25}{12}n+\frac{28}{9},
\]
and
\[
s_n^{(1)} \approx 1+\frac{3}{n}+\frac{1}{n^2}-\frac{41}{9n^3}, \quad
s_n^{(2)}\approx 1+\frac{6}{n}+\frac{6}{n^2}.
\]
Notice that for
\[
f_n^{(1)} = 1+\frac{3}{n} \quad \mbox{and} \quad g_n^{(1)} = 1+\frac{3}{n}+\frac{2}{n^2},
\]
we have
\[
\frac{(f_n^{(1)}-1)(f_{n+2}^{(1)}-1)}{(1-1/g_{n+1}^{(1)})^2} - 1
= \frac{2 (21 n^2+82 n+81)}{ (3n+5)^2 n} > 0, \quad \forall\, n \ge 1.
\]
Thus we may take $k_1=2$ and $N_1 = 1$.

Now take $k_2=0$ such that
\[
f_n^{(0)} = \frac{3}{4}n^3+\frac{3}{2}n^2+\frac{25}{12}n+\frac{19}{9},
\]
and
\[
g_n^{(0)} = \frac{3}{4}n^3+\frac{3}{2}n^2+\frac{25}{12}n+\frac{37}{9}.
\]
It can be verified that for $n \ge 8$,
\[
f_{n+1}^{(0)}/g_n^{(0)}  \ge f_n^{(1)} \quad \mbox{and} \quad g_n^{(1)} \ge g_{n+1}^{(0)}/f_n^{(0)}.
\]
Therefore, we may take $N_2=8$.

Clearly, $R_0(n)$ is negative and $R_1(n),R_2(n),R_3(n)$ is positive for each $n \ge 0$. It is routine to check that \eqref{u} and \eqref{v} hold for $n \ge 0$. Finally, for $n=5,6,7,8$ we have
\[
f_n^{(0)} \le a_{n+1}/a_n \le g_n^{(0)}.
\]
Therefore, the sequence $\{H_n(3)\}_{n \ge 8}$ is 2-log-convex. \qed

We implement a {\tt Mathematica} package to do the above computation automatically. With this package, we can reprove almost all the results in \cite{Dos08} and \cite{Chen11} automatically.

\section*{Acknowledgements}

We wish to thank professor Authur L. B. Yang for his conjecture that $\{H_n(3)\}_{n \ge 2}$ is log-convex. This work was supported by the National Science Foundation of China.

\end{document}